\documentclass[a4paper]{amsart}

\usepackage[english]{babel}
\usepackage[latin1]{inputenc}

\usepackage{eucal, longtable, amscd, amssymb, amsthm, amsmath,
  stmaryrd, calc, multicol, lscape, varioref, array}

\theoremstyle{plain}
\newtheorem*{conj}{Conjecture}

\newtheorem{thm}{Theorem}
\newtheorem*{coro}{Corollary}
\newtheorem{prop}{Proposition}
\newtheorem{prope}{Property}
\newtheorem{lemma}{Lemma}

\theoremstyle{definition}
\newtheorem*{defi}{Definition}
\theoremstyle{remark}
\newtheorem*{rem}{Remark}

\newtheorem*{eg}{Example}

\newcommand{\VSpace}{\rule{0pt}{12pt}}

\newcommand{\tablespace}{\vskip 2cm}

\DeclareMathOperator{\HH}{{\sf H}} 
\DeclareMathOperator{\Hom}{{\sf Hom}}
\DeclareMathOperator{\Pic}{{\sf Pic}}
\DeclareMathOperator{\QH}{{\sf QH}} 
\DeclareMathOperator{\RRR}{{\sf R}} 
\DeclareMathOperator{\Sym}{Sym}
\DeclareMathOperator{\vd}{{\sf virt.dim.}} 
\newcommand{\Bb}{\mathcal{B}}
\newcommand{\CC}{\mathbb{C}}
\newcommand{\Cc}{\mathcal{C}}
\newcommand{\Chow}{\mathsf{CH}}
\newcommand{\Dd}{\mathcal{D}}
\newcommand{\Ee}{\mathcal{E}}
\newcommand{\Ff}{\mathcal{F}}
\newcommand{\Gg}{\mathcal{G}}
\newcommand{\Mzeroenne}{\overline{M}_{0,n}}
\newcommand{\Mzerotre}{\overline{M}_{0,3}}
\newcommand{\Mzerozero}{\overline{M}_{0,0}}
\newcommand{\NN}{\mathbb{N}}
\newcommand{\Oo}[1][]{{\mathcal{O}_{#1}}}
\newcommand{\PP}{\mathbb{P}}
\newcommand{\Th}{\hat{T}}
\newcommand{\Xt}{\tilde{X}}
\newcommand{\ZZ}{\mathbb{Z}}
\newcommand{\fib}{\varphi}
\newcommand{\ot}{\leftarrow}
\newcommand{\pt}{\bullet}
\newcommand{\ret}{\varrho}

\title[Quantum Cohomology of some Fano Threefolds\dots]{On the quantum
  cohomology\\of some Fano threefolds\\and a conjecture of Dubrovin}

\author[G. Ciolli]{Gianni Ciolli} 

\address{Dipartimento di Matematica ``Ulisse Dini''\\Viale Morgagni 67a\\50134 Firenze, ITALY}
\email{ciolli@math.unifi.it}
\urladdr{http://www.math.unifi.it/{\textasciitilde}ciolli}


\begin{document}
\begin{abstract}
  In the present paper the small Quantum Cohomology ring of some Fano
  threefolds which are obtained as one- or two-curve blow-ups from
  $\PP^3$ or the quadric $Q^3$ is explicitely computed. Because of
  systematic usage of the associativity property of quantum product
  only a very small and enumerative subset of Gromov-Witten invariants
  is needed.  Then, for these threefolds the Dubrovin conjecture on
  the semisimplicity of Quantum Cohomology is proven by checking the
  computed Quantum Cohomology rings and by showing that a smooth Fano
  threefold $X$ with $b_3(X)=0$ admits a complete exceptional set of
  the appropriate length.
\end{abstract}

\maketitle

\begingroup
\parindent 0pt
\parskip 5pt

\section*{Introduction}

Throughout this paper, we call {\em Fano threefold} a smooth complex
compact algebraic variety of dimension three whose anticanonical
divisor is ample.

Fano threefolds have been classified by Iskovskih, Mori and Mukai
\cite{iskovskih-1977,iskovskih-1978,mori-mukai-1982,mori-mukai-2003}
in 106 deformation classes, and the second Betti number $b_2$ ranges
from $1$ to $10$.  According to the classification, we denote with
$M^b_n$ the $n$-th element of the list of Fano threefolds $M$ having
$b_2(M)=b$.

In the present paper we compute an explicit presentation for the small
Quantum Cohomology ring of 13 Fano threefolds which can be constructed
as blow-ups of $\PP^3$ or $Q^3$ along one or two smooth rational
curves. Two of them were already studied, since they can also be
constructed as projective bundles; however, we include them both for
the sake of completeness and because our computations are quite
simple.
\begin{thm}\label{thm:b2=2}
  The small Quantum Cohomology ring of the Fano threefold $M^2_k$ with
  $k=21,22,26,27,29,30,33$ is isomorphic to the polynomial quotient
  ring
  $$\CC[E,H,q_0,q_1]\Big/ (f_1^C-f_1^Q,f_2^C-f_2^Q)\ ,$$
  where the
  relations are described in table \ref{tab:rels-2}.
\end{thm}
\begin{thm}\label{thm:b2=3}
  The small Quantum Cohomology ring of the Fano threefold $M^3_n$ for
  $n=10,12,15,18,20,25$ is isomorphic to the polynomial quotient ring
  $$\CC[E_1,E_2,H,q_0,q_1,q_2]\Big/
  (f_1^C-f_1^Q,f_2^C-f_2^Q,f_3^C-f_3^Q)\ ,$$
  where the relations are
  described in table \ref{tab:rels-3}.
\end{thm}  

In section \ref{sect:exceptional} we focus on a conjecture of
Dubrovin:
\begin{conj}[Dubrovin, \cite{dubrovin-1998}~4.2.2~(1);
  \cite{bayer-manin-2001}]
  \label{conj:dubrovin}
  Let $X$ be a (smooth complex compact) variety.  The even Quantum
  Cohomology ring of $X$ is generically semisimple if and only if\\
  (i) $X$ is Fano and\\
  (ii) the bounded derived category of coherent sheaves on $X$ admits
  a complete exceptional set whose length equals
  $\sum_q\dim\HH^{q,q}(X).$
\end{conj}

Dubrovin conjecture has been proved by Bayer and Manin
\cite{bayer-manin-2001} for Del Pezzo surfaces; moreover, it holds for
$\PP^n$, because of Be\u{\i}linson's Theorem and well-known Quantum
Cohomology computations.

We prove the conjecture for other 35 Fano threefolds:
\begin{thm}\label{thm:dc}
  Dubrovin conjecture holds for the following 36 Fano threefolds (out
  of the 59 having only even cohomology):\\
  $\PP^3,Q^3,V_5,V_{22}$;\\
  $M^2_k$ with $k=21,22,24,26,27,29,30,31,32,33,34,35,36$;\\
  $M^3_k$ with $k=10,12,15,17,18,20,24,25,27,28,30,31$;\\
  $\PP^1\times S_k$ with $k=2,\dots,8$.
\end{thm}
The proof follows from Theorems~\ref{thm:b2=2} and \ref{thm:b2=3},
existing computations \cite{ancona-maggesi-2002, costa-miro-roig-2000,
  qin-ruan-1998, spielberg-1999} and the following proposition.
\begin{prop}\label{prop:exceptional}
  Let $X$ be a smooth Fano threefold with $b_3(X)=0$; then the bounded
  derived category of coherent sheaves on $X$ admits a complete
  exceptional set of length $\sum_q\dim\HH^{q,q}(X).$
\end{prop}

This method can eventually lead to a complete result for all the 59
Fano threefolds having only even cohomology; on the other hand, the
remaining 47 Fano threefolds, i.~e. those $X$ having $b_3>0$, present
additional difficulties.

Precisely, many of such threefolds are known not to have a complete
exceptional set of the appropriate length; thus, the Dubrovin
conjecture is equivalent to non-semisimplicity of the big Quantum
Cohomology ring. The computation of this latter ring is much more
complicated, since it involves infinitely many Gromov-Witten
invariants; also, the relation between the small ring and
non-semisimplicity of the big ring is rather implicit, contrarily to
the case of Fano threefolds having $b_3=0$.

\subsubsection*{Acknowledgement} I wish to thank V.~Ancona for many
useful discussions and for having introduced me into the present
subject. I am also grateful to the {\em Dottorato di Ricerca} of the
University of Florence and to the Italian MIUR Project {\em
  Propriet\`a geometriche delle variet\`a reali e complesse} for their
financial support during the preparation of my {\em Tesi di
  Dottorato}, when some of the results in the present paper were
obtained. Moreover, I express my gratitude to Tom Bridgeland for
having noticed a small mistake on a previous version of this paper.

\section{Quantum Cohomology}

Let $X$ be a smooth projective $n$-dimensional variety, and let
$p=b_2(X)$. A {\em Quantum Cohomology ring of $X$} is a tangent space
to the classical cohomology ring $\HH^*(X,\CC)$ endowed with a
so-called {\em quantum} multiplicative structure
\cite{kontsevich-manin-1994} which depends on the tangency point
$\gamma\in\HH^*(X,\CC)$.

All the Quantum Cohomology rings with $\gamma\in\HH^2(X,\CC)$ can be
described simultaneously by a parametric presentation, which depends
on the $p$ parameters $q_1,\dots,q_p\in\CC\setminus\{0\}$ and gives
the so-called {\em small Quantum Cohomology ring}. In this paper we
are concerned only with the small Quantum Cohomology ring of $X$,
which we will denote by $\QH(X)$; for simplicity, in the sequel we
will drop the word ``small''.

Let $T_0=1=[X], T_1,\dots,T_m$ be a homogeneous basis of the graded
vector space $\HH^*(X,\CC)$ such that $|T_i|\leq|T_j|$ if $i<j$. We
denote by $\pt$ the cohomology class of a point in $X$. Let
$\Th_0=\pt,\dots,\Th_m\in\HH^*(X,\CC)$ be the basis which is
Poincar\'e dual to the $T_i$, that is, the basis such that $(T_i\cup
\Th_j)_n=\delta_{ij}\pt$ for all $i,j$. Let $\Bb\subset\HH_2(X,\ZZ)$
denote the subset of nonzero homology classes of effective rational
maps to $X$. Fix $\beta_1,\dots,\beta_p\in \HH_2(X,\ZZ)$ forming a
vector space basis of $\HH_2(X,\CC)$ and such that any $\beta\in\Bb$
can be expressed as $\sum_i\beta_i b_i$ with $b_i\in\NN$.

\begin{defi}
  Let $n\geq3$. The {\em moduli space $\Mzeroenne(X,\beta)$ of
    genus-zero $n$-marked maps of class $\beta$ in $X$} is the moduli
  space of isomorphism classes of stable $n$-pointed maps
  $(f,\Cc,x_1,\dots,x_n)$ such that $f_*[\Cc]=\beta$.
\end{defi}
Each space $\Mzeroenne(X,\beta)$ is endowed with $n$ natural {\em
  evaluation maps} to $X$
$$i_k:(f,\Cc,x_1,\dots,x_n)\mapsto f(x_k)\ $$
which are involved in
the following
\begin{defi}
  The {\em Gromov-Witten invariant} $I_\beta(\gamma_1,\dots,\gamma_n)$
  is the rational number corresponding to the value of the degree-zero
  component of the cap product
  $$Z\cap(i_1^*\gamma_1\cup\dots\cup i_n^*\gamma_n)$$
  where
  $Z\in\Chow_*(\Mzeroenne(X,\beta))$ is the {\em virtual fundamental
    class}.
\end{defi}

Gromov-Witten invariants satisfy some properties
\cite{kontsevich-manin-1994} among which we recall the following two.
\begin{prope}[Divisor Axiom]
  Let $D$ be a divisor; then
  $$I_\beta(D,\alpha_1,\dots,\alpha_n)=(\beta.D)I_\beta(\alpha_1,\dots,\alpha_n)\
  .$$
\end{prope}
\begin{prope}[Grading Axiom]
  Let $\alpha_1,\dots,\alpha_n$ be homogeneous cohomology classes,
  i.~e. $\alpha_i\in\HH^{\deg\alpha_i}(M,\CC)$. If
  $$\vd(\Mzeroenne(M,\beta))\neq \sum_i\deg\alpha_i$$
  then the
  Gromov-Witten invariant $I_\beta(\alpha_1,\dots,\alpha_n)$ vanishes.
\end{prope}

\begin{rem}[$n$-pointed invariants with $n<3$]
  We can define Gromov-Witten invariants even in the case $n=1,2$ by
  means of the Divisor axiom: we define respectively
  $I_\beta(\gamma_1)$ and $I_\beta(\gamma_1,\gamma_2)$ as the rational
  numbers $x_1,x_2$ such that
  $I_\beta(\gamma_1,D_2,D_3)=x_1(D_2.\beta)(D_3.\beta)$ and
  $I_\beta(\gamma_1,\gamma_2,D_3)=x_2(D_3.\beta)$ where $D_1,D_2$ are
  divisors such that $\beta.D_i\neq 0$. This definition is easily seen
  to be independent on the choice of the $D_i$.
\end{rem}

Given $\alpha_1,\alpha_2\in\HH^*(X,\CC)$, their quantum product is
defined as
$$\alpha_1*\alpha_2=\alpha_1\cup\alpha_2+\sum_{\beta\in\Bb} \sum_\ell
I_\beta(\alpha_1,\alpha_2,T_\ell)\Th_\ell q^\beta$$
where
$q^\beta:=q_1^{b_1}\cdots q_p^{b_p}$ and the
$I_\beta(\gamma_1,\gamma_2,\gamma_3)$ are the three-pointed
Gromov-Witten invariants, which encode enumerative information about
the geometry of rational curves on $X$.

This product extends $\CC[q_1,\dots,q_p]$-linearly to the ring
$\HH^*(X,\CC)[q_1,\dots,q_p]$. It is widely known that if $X$ is Fano,
as we assume in the present paper, the sum $\sum_\beta$ is finite
since almost all the $I_\beta$ vanish because of the Grading axiom
(see below).

Quantum product has been proven to be associative
\cite{ruan-tian-1995,behrend-1997,behrend-fantechi-1997,li-tian-1998};
this property is equivalent to a large family of relations between
Gromov-Witten invariants, which we use later to reduce the amount of
the enumerative data which is essential for the determination of
quantum multiplication.

These properties lead us to the following definition.
\begin{defi}
  We will call a Gromov-Witten invariant
  $I_\beta(\alpha_1,\dots,\alpha_n)$ {\em essential} if
  \begin{itemize}
  \item $\deg\alpha_i>2$ for all $i$, that is, Divisor Axiom cannot be
    applied any further;
  \item $\vd(\Mzeroenne(M,\beta))= 2\sum\deg\alpha_i$, that is, it
    does not vanish trivially because of Grading Axiom;
  \item it is involved by some associativity relation in the small
    Quantum Cohomology ring, that is, it appears as a coefficient in
    some expression of type $T_i*(T_j*T_k)-(T_i*T_j)*T_k$ with $\deg
    T_i+\deg T_j+\deg T_k\leq 2\dim X$.
  \end{itemize}
\end{defi}

From the definition it follows that essential invariants determine
also non-essential ones. Thus the knowledge of all essential
Gromov-Witten invariants completely determines Quantum Cohomology.
Indeed, as showed in the following proposition, from the knowledge of
quantum multiplication we can extract a presentation for the ring
$\QH(X)$ as a quotient of the polynomial ring
$\HH^*(X,\CC)[q_1,\dots,q_p]$.

\begin{prop}[\cite{fulton-pandharipande-1997}, Proposition 11]\label{prop:qh-ring}
  Let $f_1,\dots,f_r\in\CC[x_1,\dots,x_k]$ be homogeneous polynomials
  such that $$\HH^*(X,\CC)\cong\CC[x_1,\dots,x_k]\ \big/\ 
  \Big(f_1,\dots,f_r\Big)$$
  as graded rings. Let $f_1',\dots,f_r'$ be
  homogeneous elements in $\CC[x_1,\dots,x_k,q_1,\dots,q_p]$, where
  for all $i$ $\deg q_i=-K_X.\beta_i$ and
  \begin{itemize}
  \item $f_i'(x_1,\dots,x_k,0,\dots,0)=f_i(x_1,\dots,x_k)$\ ;
  \item $f_i'(x_1,\dots,x_k,q_1,\dots,q_p)=0$ in
    $\HH^*(X,\CC)\otimes\CC[q_1,\dots,q_p]$ with the quantum product.
  \end{itemize}  
  Then
  $$\QH^*(X,\CC)\cong\CC[x_1,\dots,x_k,q_1,\dots,q_p]\ \big/\ 
  \Big(f_1',\dots,f_r'\Big)\ .$$
\end{prop}

\section{Existing computations}
A presentation of the Quantum Cohomology ring has been computed
\cite{qin-ruan-1998, spielberg-1999, costa-miro-roig-2000,
  ancona-maggesi-2002} for the 21 Fano threefolds which can be
constructed as $\PP(\Ee)$ \cite{szurek-wisniewski-1990}, where $\Ee$
is a rank-2 vector bundle over a surface. Also, 18 Fano threefolds
have a toric representative \cite{batyrev-1981,watanabe-watanabe-1982}
and their Quantum Cohomology can be studied also with the techniques
available for toric varieties
\cite{batyrev-1993,givental-1998,spielberg-1999}. These two classes
overlap; their union counts 28 Fano threefolds. 

Moreover, Quantum Cohomology has been computed \cite{beauville-1995}
for some complete intersections in projective spaces, under some
hypotheses which hold for three Fano threefolds: the quadric $M^1_2$,
the cubic $M^1_5$ and $V_5=M^1_7$.

A few of our computations were known before. The two threefolds
$M^2_{33}$ and $M^3_{25}$ are toric and could have been treated with
the above-mentioned techniques.  Quantum Cohomology of the threefold
$M^2_{27}$ was already computed in \cite{ancona-maggesi-2002} since
$M^2_{27}$ can also be constructed as a projective bundle over a
surface.  However, we report such presentations in order to underline
the analogy with the other blow-ups.

Some Gromov-Witten invariants of blow-ups along smooth curves (in
particular, lemmas~\ref{lemma:L0:2:21,26,29}
and~\ref{lemma:L0:3:10,15,20}) follow also from Theorem 1.5 in
\cite{hu-2000}, in which symplectic techniques are used to prove that
$$I^X_\beta(\alpha_1,\dots,\alpha_m) =
I^{\Xt}_{f!(\beta)}(f^*\alpha_1,\dots,f^*\alpha_m)$$
where
$f!(\beta):=Pf^*P\beta$, $P$ is Poincaré duality and
$\beta\in\HH_2(X,\ZZ)$; anyway, we present our enumerative proofs,
which happen to be particularly simple and clear.

We check Dubrovin conjecture by means of explicit presentations of
Quantum Cohomology rings; but for some manifolds only Gromov-Witten
invariants or a sketch of the computations have been done before. This
is the case for the threefolds $Q^3,V_5,V_{22}$ and for trivial
$\PP^1$-bundles over Del Pezzo surfaces, for which we computed
explicit Quantum Cohomology presentations starting from partial
results in \cite{ancona-maggesi-2002, bayer-manin-2001}.

Very recently, A.~Bayer proved in \cite{bayer-2004} that
semisimplicity of Quantum Cohomology is preserved by point blow-ups,
and proposed consequently a small change to the conjecture, that is,
to remove the condition of being Fano.

\section{The blow-up of a Fano threefold}
We recall the following description of the even cohomology ring of the
blow-up of a threefold along a smooth curve (see e.g. lemma 2.11 in
\cite{iskovskih-1977}, or \cite{fulton-1984}).

\begin{lemma} \label{lemma:chowblow}
  Let $f:\Xt\to X$ be the blow-up of a smooth threefold along a smooth
  curve $\Cc$. Let $E\in \HH^2(\Xt,\CC)$ be the class of the
  exceptional divisor, and let $\fib\in \HH^4(\Xt,\CC)$ be the class
  of an exceptional fiber.  Then
  $$\HH^*(\Xt,\CC)\cong f^*\HH^*(X,\CC)\oplus \CC E\oplus \CC\fib$$
  as
  vector spaces, with the multiplicative structure defined by
  \begin{itemize}
  \item $E^2=-f^*[\Cc]+c_1\fib$,
  \item $E.\fib=-1$,
  \item $E.f^*D=(\Cc.D)\fib$ and $\fib.f^*D=0$ for all $D\in\HH^2(X)$,
  \item $E.f^*C=\fib.f^*C=0$ for all $C\in\HH^4(X)$,
  \end{itemize}
  where $c_1:=\deg c_1(N_{\Cc}X)=2g(\Cc)-2-K_X.[\Cc]$.
\end{lemma}

Also the Quantum Cohomology structure of $\Xt$ is closely related to
the one in $X$. Indeed, an irreducible rational curve in $\Xt$ can be
either the strict transform of a rational curve in $X$ or an
exceptional fiber; starting from this argument, we can build a family
of finite morphisms between the various components of moduli spaces of
stable maps of $X$ and $\Xt$, expressing all the genus-zero
Gromov-Witten invariants of $\Xt$ by means of the ones of $X$, as in
the example \vpageref{eg:m-2-30}.

Moreover, a blow-up of the above type always gives a ``standard''
Gromov-Witten invariant, as stated in the next lemma; this invariant
turns out to be always essential.

\begin{lemma}\label{lemma:I_F(fib)}
  In a threefold $\Xt$ which is the blow-up along a curve, the class
  $F$ of an exceptional fiber is enumerative and the value of the
  Gromov-Witten invariant $I_F(\fib)$ is $-1$, where $\fib$ is the
  cohomology class of an exceptional fiber.
\end{lemma}
\begin{proof}
  Exceptional fibers are parameterized by the blown-up curve, and a
  curve in the class $F$ must be an exceptional fiber. So the
  dimension of the moduli space $\Mzerotre(X,F)$ is
  $$1+3=4=3+\Big(\sigma^*(-K)-E\Big).F=3+0+1\ .$$
  
  Let $\ret_i$ be the cohomology class of the strict trasform of a
  line which intersects the blown-up curve transversally in $i$
  points; we have $I_F(\ret_i)=i$; since $\ret_0=\ret_1+\fib$, the
  claim follows from (multi)linearity of Gromov-Witten invariants.
\end{proof}

\section{Getting information from associativity}\label{sect:associativity}
Only finitely many Gromov-Witten invariants are involved in the
computation of all the quantum products needed in Proposition
\ref{prop:qh-ring}.

Divisor and Grading axioms (see \cite{kontsevich-manin-1994}) together
with the fact that $b_1(X)=0$ (true for any Fano threefold) reduce the
family of ``essential'' invariants to those
$I_\beta(\alpha_1,\dots,\alpha_n)$ such that $n\leq 3$,
$\deg\alpha_i>2$ for all $i$ and $\sum\deg\alpha_i=-K.\beta+n$\ .
Let's label such invariants as $x_1,\dots,x_N$.

Consider a threefold with $b_1=0$, e.~g. a Fano threefold; then the
only non-trivial cup products are those between three divisors. This
is not true for quantum product; indeed, the small Quantum Cohomology
ring has infinite dimension if regarded as a $\CC$-vector space.
However, to reduce the size of the computations we considered only
those associativity relations arising from triple products of
divisors; they involve only a subclass of Gromov-Witten invariants,
which we hope to be enough large to fulfill our purpose. This is in
the spirit of the Reconstruction theorem (3.1 of
\cite{kontsevich-manin-1994}), which implies that the complete family
of Gromov-Witten invariants is determined by the above subclass.

Consider all the possible triples $1\leq i\leq j\leq k\leq p=b_2(X)$
where at least two of $i,j,k$ are distinct, and write
\begin{equation}
  \label{eq:Pijk}
  P(i,j,k):=(T_i*T_j)*T_k-T_i*(T_j*T_k)=0\ .
\end{equation}

$P$ is a polynomial in $\CC[T_1,\dots,T_p,x_1,\dots,x_N]$; if we
decompose it as $P(i,j,k)=\sum_I P(i,j,k)_I T^I$, where
$I=(i_1,\dots,i_p)$ is a multi-index, $T^I=\prod_s T_s^{i_s}$ and
$P(i,j,k)_I\in\CC[x_1,\dots,x_N]$ we obtain a system of associativity
equations
\begin{equation}
  \label{eq:assoc}
  \begin{cases}
    P(i,j,k)_I=0\\
    \text{for all multi-indices $I$}\\
    \text{for all non-diagonal triples }(i,j,k)\in\Sym^3\{1,\dots,p\}\ ;
  \end{cases}
\end{equation}
the ideal $J_A$ generated by such polynomials contains all the
possible associativity relations among the $x_i$.

In general these relations do not suffice to determine all the
invariants, since the affine variety $A\subset\CC^N$ cut out by
system~(\ref{eq:assoc}) has positive dimension; however, by
intersecting $A$ with a suitable affine variety $G$ defined by some
additional geometric relations
$J_G=(g_1,\dots,g_s)\subset\CC[x_1,\dots,x_N]$, we can have that
$G\cap A=\{\text{one point}\}$, that is, we can determine the Quantum
Cohomology of $X$ by using the minimal geometric information expressed
by $G$.

To carry out this procedure, we used a specific computer program
\cite{algheme,fano3fold-software} which computes all the equations in
\eqref{eq:assoc} and builds a list of all the essential invariants.
Our purpose is to determine $J_G$, that is, a set of geometric
relations which suffice to determine all the Gromov-Witten invariants
of the manifold. 

We begin with an empty ideal $J_G=(0)$. We compute a {\em geometric
  relation} which we don't have already, i.~e. which is not contained
in $J_G+J_A$ (e.~g. we compute the value of some Gromov-Witten
invariant not determined by associativity); we add this relation to
$J_G$ and we check whether $J_G+J_A$ is the ideal of a single point in
$\CC^N$. If not, we repeat this procedure until $J_G+J_A$ is as big as
desired.

In the thirteen cases that we studied (cfr.
Table~\ref{tab:geom-info}), this procedure yielded the desired results
after a few steps; moreover, the standard relations coming from the
blow-up construction (lemma~\ref{lemma:I_F(fib)}) can be used as
geometric relations. We refer to the following example and to
\cite{ciolli-phd-thesis} for further details about the computational
part.

\begin{eg}\label{eg:m-2-30}
  Consider the blow-up $M=M^2_{30}$ of $\PP^3$ along a conic $\Cc$.
  Classical relations are $E^2-3EH+2H^2=0$ and $EH^2$; to quantize
  them we only need classes $\beta$ such that $3\leq\vd\beta\leq6$,
  that is, such that $3\geq -K_M.\beta=(4H-E).(dL_0+f F)=4d+f$, where
  $L_0=f^*[\text{line}]$ and $F$ is the class of an exceptional fiber.
  
  Let $f(d)$ be the minimum integer $f$ such that there exists a
  rational curve $Z\subset\PP^3$ whose homology class is
  $[Z']=dL_0+fF$. Since $f(d)=-2d$, we put $L=L_0-2F$; the only
  $\beta$'s satisfying $-K_M.\beta\leq 3$ are $F$, $2F$, $3F$, $L$ and
  $L+F$.
  
  ``Exceptional'' moduli spaces $\Mzerozero(M,kF)$ are branched
  coverings of $\Cc$, and are easily seen to have expected dimension,
  as the other two. Indeed, $\Mzerozero(M,L)$ is isomorphic to the
  dual of the plane spanned by $\Cc$, and $\Mzerozero(M,L+F)$ has two
  components: the two-dimensional one is composed by reducible curves,
  and is a $2:1$ covering of $\Mzerozero(M,L)$, while irreducible
  curves lie in the three-dimensional component $Z$, because the map
  $\big\{(x,y)\in\PP^3\times\Cc\big|x\neq y\big\} \to Z$ %
  sending $(x,y)$ to the strict transform of the line $<x,y>$ is a
  fibration with one-dimensional fibers.
  
  Associativity equations involve $N=14$ Gromov-Witten invariants; the
  affine variety $A$ has dimension three and meets in a single point
  the variety $G$ defined by $I_F(\fib)=-1$, $I_L(\fib,\fib)=1$ and
  $I_{F+L}(\fib,\pt)=2$.  The first invariant is standard, while the
  other two are easily computed considering the above description of
  the related moduli spaces.  So Quantum Cohomology of $M$ is
  completely determined by associativity constraints and ``simple''
  enumerative geometric information.
\end{eg}

\section{Exceptional objects in the Derived Category\\and the Dubrovin
  Conjecture}\label{sect:exceptional}

We recall some terminology and results from \cite{orlov-1993}.

\begin{defi}
  Let $X$ be a smooth complex projective variety. We denote with
  $\Dd(X)$ the bounded derived category of coherent sheaves on
  $X$. 
\end{defi}

\begin{defi}
  An object $E$ in $\Dd(X)$ is said to be {\em exceptional} if
  $$\RRR^i\Hom(E,E)=
  \begin{cases}
    0   & \text{if $i>0$}\ ,\\
    \CC & \text{if $i=0$}\ .
  \end{cases}$$
\end{defi}

\begin{defi}
  An {\em exceptional set} (also called a {\em system of exceptional
    objects} in $\Dd(X)$ of length $n+1$ is an ordered set
  $$(E_0,\cdots,E_n)$$
  such that
  \begin{itemize}
  \item all the $E_i$ are exceptional;
  \item $\RRR^*\Hom(E_i,E_j)=0$ if $i>j$ (the {\em semiorthogonality}
    condition).
  \end{itemize}
  Moreover, the set (system) is called {\em complete} (full) if it
  generates $\Dd(X)$ as a triangulated category.
\end{defi}

Be\u{\i}linson's theorem exhibits a complete exceptional set for the
projective space $\PP^n$.
\begin{thm}[Be\u{\i}linson
  \cite{beilinson-1978,ancona-ottaviani-1989}]\label{thm:beilinson}
  Given a coherent sheaf $\Ff$ on $\PP^n$, there exists a finite
  complex of sheaves
  $$\Gg^\bullet:\qquad0\to\Gg^{-n}\to\dots\to\Gg^n\to 0$$
  such that
  each one of the $\Gg^q$ is the direct sum of bundles of type
  $\Oo(-n),\dots,\Oo$ and
  $$\HH^q(\Gg^\bullet)=\begin{cases}
    0   & \text{if $q\neq 0$,}\\
    \Ff & \text{if $q=0$.}
  \end{cases}$$
\end{thm}

Indeed, in the derived category any object is equivalent to anyone of
its resolutions, so these $n+1$ sheaves generate $\Dd(\PP^n)$. Also,
the exceptionality of the $\Oo(i)$ is an easy consequence of the fact
that $\RRR^0F=F$ for any functor $F$ and that the $\Oo(i)$ are free,
so that $\RRR^q\Hom(\Oo(i),\Oo(j))=0$ for all $q>0$ and all $i,j$.
These arguments also imply semiorthogonality of the set.

Dubrovin conjecture holds for $\PP^n$; indeed,
$\QH(\PP^n)\cong\CC[H,q]/(H^{n+1}-q)$ and a corollary of
Be\u{\i}linson's theorem is that the length of the exceptional set is
exactly
$$\dim\HH^*(\PP^n,\CC)=\sum_q\dim\HH^{q,q}(\PP^n)=n+1\ .$$
Conforming to
Dubrovin conjecture, we will say that the {\em appropriate} length for
an exceptional set of a manifold $X$ is $\sum_q\dim\HH^{q,q}(X)$.

Starting from manifolds whose derived category has a complete
exceptional set of the appropriate length, as is $\PP^n$, we can
construct other manifolds with the some property if we use processes
which preserve both the existence of an exceptional set and the
appropriatedness of its length.  This is the case for some blow-ups
and for all projectivizations of vector bundles, as we can deduce from
the following two theorems.

\begin{thm}[\cite{ancona-ottaviani-1991}; Cor.~2.7 of
  \cite{orlov-1993}] \label{thm:good-proj} Let $\Ee\to M$ be a vector
  bundle over a manifold.  If $\Dd(M)$ has a complete exceptional set,
  then also $\Dd(\PP(\Ee))$ possesses a complete exceptional set.
\end{thm}

\begin{coro}
  If $\Dd(M)$ has a complete exceptional set of the appropriate
  length, then also $\Dd(\PP(\Ee))$ has a complete exceptional set of
  the appropriate length.
\end{coro}
\begin{proof}
  The statement about the length is a straightforward computation that
  can be obtained from one of the proofs of the theorem in
  \cite{ancona-ottaviani-1991,orlov-1993} and from Leray-Hirsch
  theorem describing the cohomology of a projective bundle.
\end{proof}

\begin{thm}[Cor.~4.4 of \cite{orlov-1993}] 
  \label{thm:good-blowup}
  Let $M$ be the blow-up of a smooth variety $X$ along a smooth
  subvariety $Y$.  If $\Dd(X)$ and $\Dd(Y)$ have a complete
  exceptional set, then also $\Dd(M)$ possesses a complete exceptional
  set.
\end{thm}

\begin{coro}
  Let $M$ be the blow-up of a threefold $X$ along a smooth irreducible
  subvariety $Y$; suppose that $\Dd(X)$ has a complete exceptional set
  of the appropriate length, and that $Y$ is either a point or a
  rational curve. Then also $\Dd(M)$ has a complete exceptional set of
  the appropriate length.
\end{coro}
\begin{proof}
  Again, the length of the exceptional set can be computed easily from
  the contents of section 4 in \cite{orlov-1993}, while the
  appropriatedness of the length follows from lemma
  \ref{lemma:chowblow}.
\end{proof}

\begin{proof}[Proof of Proposition~\ref{prop:exceptional}]  
  The classification of smooth Fano threefolds states that Fano
  threefolds having $b_3=0$ are $4+21+34=59$ out of 106. Indeed, they
  can be subdivided in the following three classes:
  \begin{itemize}
  \item the four threefolds with $b_2=1$, that is, $\PP^3,Q^3,V_5$ and
    $V_{22}$;
  \item the twenty-one threefolds which are $\PP^1$-bundles over
    surfaces;
  \item the thirty-four threefolds which are not $\PP^1$-bundles over
    surfaces and can be obtained by blowing up another Fano threefold
    along a smooth subvariety $Y$ having only cohomology of even
    degree.
  \end{itemize}
  
  Fano threefolds with $b_2=1$ for which (cfr.
  \cite{bayer-manin-2001}) a complete exceptional set in $\Dd(X)$ is
  known to exist are $\PP^3$, $Q^3$, $V_5$ \cite{orlov-1991} and
  $V_{22}$ \cite{kuznetsov-1996,faenzi-2003-v22}; they are the only
  ones with $\Pic=\ZZ$ having only even cohomology. Moreover, the
  exceptional set has the appropriate length.
  
  All the 21 Fano threefolds $M_1,\dots,M_{21}$ which are
  $\PP^1$-bundles over surfaces, corresponding to the 21 Fano vector
  bundles in \cite{szurek-wisniewski-1990} (we refer to that paper or
  to \cite{ancona-maggesi-2002} for the notation), have a complete
  exceptional set in $\Dd(M_i)$ of the appropriate length.
  
  Indeed, for all the surfaces $S$ which appear as bases of such
  bundles $\Dd(S)$ has a complete exceptional set of the appropriate
  length, since they are $\PP^2$, $\PP^1\times\PP^1$ -- which can be
  thought as the projectivization of trivial vector bundles over
  $\PP^1$ -- and Del Pezzo surfaces, for whom a computation quite
  similar to the one in Theorem~\ref{thm:good-blowup} can be made.
  Thus only a further application of the corollary of
  Theorem~\ref{thm:good-proj} proves the claim.
  
  Finally, the corollary of Theorem~\ref{thm:good-blowup} shows
  inductively the existence of a complete exceptional set of the
  appropriate length in $\Dd(X)$ for any $X$ in the third class.
\end{proof}

\section{Proofs of the main results}

We denote by $\pt$, $\fib$ and $\ret$ the cohomology class
respectively of a point, of an exceptional fiber and of the pullback
of a generic line class in $\PP^3$ or $Q^3$; we put $L_0=(\ret)_*$ and
$F=(\fib)_*$, where $(-)_*:\HH^p(M,\CC)\to\HH_p(M,\CC)$ denotes
Poincar\'e duality.

\begin{proof}[Proof of Theorem~\ref{thm:b2=2}]
  Lemma \ref{lemma:chowblow} gives the classical relations.
    
  The maximal $G$, as in section~\ref{sect:associativity}, is defined
  by the {\em standard} equation $I_F(\fib)=-1$ (given by lemma
  \ref{lemma:I_F(fib)}) plus the following equations, for which we
  refer to the lemmas in the last section:
  \begin{align}
    \label{eq:b2:L0-F}
    &I_{L_0-F}(\ret,\pt)=\deg\Cc & \text{for $n=22,27,30,33$;}\\
    \label{eq:b2:L0-2F}
    &I_{L_0-2F}(\ret,\ret)=1 & \text{for $n=30$;}\\
    \label{eq:b2:L0}
    &I_{L_0}(\ret,\pt)=1\ \text{and}\ I_{L_0-F}(\ret,\fib)=1
    & \text{for $n=21,26,29$}.
  \end{align}
\end{proof}

\begin{proof}[Proof of Theorem~\ref{thm:b2=3}]
  As for the previous theorem, Lemma \ref{lemma:chowblow} gives the
  classical relations. Moreover, in these cases quadratic relations
  suffice to generate the cohomology ring, as it can be easily
  verified. Since all of these threefold are obtained by blowing up
  two disjoint curves, we have two distinct homology classes $F_1,F_2$
  of exceptional fibers and consequently two standard Gromov-Witten
  invariants of the type described in lemma \ref{lemma:I_F(fib)},
  which for all these varieties turn out to be independent: i.~e.,
  they cut out a 2-codimensional subvariety of $A$.
  
  The maximal $G$ (cfr. section~\ref{sect:associativity}) is defined by the {\em standard}
  equations $I_{F_1}(\fib_1)=I_{F_2}(\fib_2)=-1$ plus the following
  ones, for which we refer to the lemmas in the last section:
  \begin{align}
    &I_{L_0}(\ret,\pt)=1 & \text{for $n=10,15,20$},\\
    &I_{L_0-F_1-F_2}(\pt)=\deg\Cc_2 & \text{for $n=12,18,25$},
  \end{align}
  where $\Cc_2$ is the smooth rational curve such that the threefold
  is obtained by blowing up $\PP^3$ along a disjoint union
  $\text{(line)}\sqcup\Cc_2$.
\end{proof}

\begin{proof}[Proof of Theorem~\ref{thm:dc}]
  The existence of an exceptional set is positively answered for all
  these varieties by Proposition~\ref{prop:exceptional}. Generic
  semisimplicity is verified with a standard commutative algebra
  software using the explicit presentations worked out in the present
  paper and in \cite{ancona-maggesi-2002, costa-miro-roig-2000,
    qin-ruan-1998, spielberg-1999}, after the following lemma.

  \begin{lemma}
    Generic semisimplicity of Quantum Cohomology is a consequence of
    the existence of a single semisimple point
    $\gamma\in\HH^2(X,\CC)$.  This in turn is equivalent to generic
    semisimplicity over the $\HH^2(X,\CC)$.
  \end{lemma}

  To see this, consider the following diagram:
  $$\CC^{c+p}\supset Z=V(f_1',\dots,f_k')\ 
  \stackrel{\pi}{\to} %
  \ (\CC^*)^p\ %
  \stackrel{\psi}{\ot} %
  \ \HH^2(X,\CC)\ ,$$
  where $\pi$ is the projection forgetting the
  $q_i$'s and
  $\psi(\gamma_1,\dots,\gamma_p)=\big(\exp(\gamma_1),\dots,\exp(\gamma_p)\big)$.
  Semisimplicity of the small Quantum Cohomology ring is equivalent to
  the fact that $\psi(S)$ is a nonempty Zariski open subset in
  $(\CC^*)^p$, where $S$ is the locus of the points
  $\gamma\in\HH^2(X,\CC)$ such that the Quantum Cohomology
  $\QH_\gamma(X)$ is semisimple.

  Indeed, semisimple points $\gamma\in S$ correspond exactly to the
  points $q=\psi(\gamma)\in\big(\CC^*\big)^p$ such that $\pi^{-1}(q)$
  is a reduced zero-dimensional scheme.  Moreover, the fiber
  $\pi^{-1}(\gamma)$ is obtained by putting $q_i:=\psi(\gamma_i)$
  $\forall\ i$, that is, by intersecting $Z$ with a linear variety of
  codimension $p$. So $\psi(S)$ is exactly the locus in $(\CC^*)^p$
  over which the fibers of $\pi$ are reduced, which is an open Zariski
  subset. If it is empty, then all the fibers are nonreduced, that is,
  $Z$ itself is not reduced and the small Quantum Cohomology is not
  semisimple.    
\end{proof}

\begin{lemma}\label{lemma:m-2-22,27,30,33}
  $I_{L_0-F}(\ret,\pt)=\deg\Cc$ for $X=M^2_{k}$, $k=22,27,30,33$.
\end{lemma}

\begin{proof}
  Let $\ell\in\ret$ and $x\in\pt$ be generic representatives.
  
  Irreducible curves $\Cc'$ of class $L_0-F$ are strict transforms of
  lines intersecting $\Cc$ in exactly one point. The map sending
  $\Cc'\in\Mzerozero(X,L_0-F)$ to $\Cc'\cap\Cc$ is a fibration over
  $\Cc$ with bidimensional fibers.
  
  If $\deg\Cc>1$, the class $L_0-F$ contains also reducible curves; we
  want to show that these curves do not contribute to the invariant,
  i.~e., that no line intersecting $\Cc$ with multiplicity greater
  than 1 can meet both $x$ and $\ell$. 
  
  Indeed, choose a generic plane $\PP^2$ containing $\ell$.  The
  projection $\pi_x:\PP^3\to\PP^2$ maps $\Cc$ into a plane rational
  curve $\Cc_0$ with $m=\frac12(d-1)(d-2)$ nodes $x_1,\dots,x_m$. The
  lines $\overline{x_i\,x}$ are exactly all the lines passing through
  $x$ and $\ell$ and intersecting $\Cc$ with multiplicity greater than
  $1$.  Genericity of $x$ implies that none of these $m$ lines meets
  $\ell$, so that no reducible curve in class $L_0-F$ is involved in
  the computation of this Gromov-Witten invariant.
  
  This implies both the enumerativeness of the invariant, since we
  have seen above that $\dim\Mzerozero(X,L_0-F)=3$, and the
  possibility to compute it considering only irreducible curves, that
  is, counting the number of strict transforms of lines $\ell$
  touching $\Cc$ in exactly one point.
  
  The cone projecting $\Cc$ from $x$ intersects $\ell$ in $d=\deg\Cc$
  points, which belong to the lines $\ell_i$ ($i=1,\dots,d$) in the
  the cone. The strict transforms of the $\ell_i$ are thus all the
  rational curves belonging to class $L_0-F$ and meeting the strict
  transforms of $x$ and $\ell$, thus proving the claim.
\end{proof}

\begin{lemma}\label{lemma:m-2-30}
  $I_{L_0-2F}(\ret,\ret)=1$ for $M^2_{30}$.
\end{lemma}

\begin{proof}
  Lines who are bisecant a given conic in $\PP^n$ are parameterized by
  the dual of the plane spanned by the conic; thus the dimension of
  $\Mzerozero(M^2_{30},L_0-2F)$ is $2+3=5=\vd(L_0-2F)$ and we can use
  the enumerative interpretation.
  
  Let $Z$ be the plane spanned by the conic $\Cc$, and let
  $\ell_1,\ell_2$ be two generic lines in $\PP^3$; let
  $Z\cap\ell_i=\{x_1\}$.
  
  The value of $I_{L_0-2F}(\ret,\ret)$ is one, since it counts the
  number of lines $\ell$ included in $Z$ whose strict transform
  $\ell'$ meets the strict transforms of $x_1$ and $x_2$, that is, the
  number of lines in a plane passing through two generic points.
\end{proof}

\begin{lemma}\label{lemma:L0:2:21,26,29}
  $I_{L_0}(\ret,\pt)=1$ for $M^2_k$, $k=21,26,29$.
\end{lemma}

\begin{proof}
  Let $M$ be the blow-up of $Q^3$ along the rational curve $\Cc$.
  
  The family of lines in $Q^3$ is three-dimensional; the class of the
  strict transform of a line $\ell$ is $L_0$ if and only if $\ell$
  does not meet $\Cc$; it is easy to see that the subspace of
  reducible curves in $\Mzerotre(M,L_0)$ has positive codimension.

  Thus the dimension of the whole space of rational curves having
  class $L_0$ is $3+3=6=\vd(L_0)$ and the enumerative interpretation
  is legitimate, that is, the invariant $I_{L_0}(\ret,\pt)$ counts the
  number of rational curves of class $L_0$ meeting the strict
  transforms of a generic line and a generic point in $Q^3$.
  
  Consider a generic line $\ell'\subset Q^3$ and a generic point
  $x'\in Q^3$; there exists exactly one line $\ell\subset Q^3$ meeting
  both $x'$ and $\ell'$; moreover, the genericity of $x'$ and $\ell'$
  implies that $\ell$ does not meet $\Cc$. Since any curve of class
  $L_0$ is the strict transform of a line $\ell$ in $Q^3$ (plus some
  exceptional fiber if $\ell$ meets $\Cc$) $\ell$ is the only curve of
  class $L_0$ meeting both $x'$ and $\ell'$, and
  $I_{L_0}(\ret,\pt)=1$.
\end{proof}

\begin{lemma}\label{lemma:N11}
  $I_{L_0-F}(\ret,\fib)=1$ for $M^2_k$, $k=21,26,29$.
\end{lemma}
\begin{proof}
  A curve of class $L_0-F$ which meets a generic line $\ret'$ cannot
  be decomposed as $\Cc'\sqcup \fib'$, since in that case $\Cc'$ would
  be bisecant to $\Cc$ and thus could not meet $\ret'$.
  
  An irreducible curve of class $L_0-F$ is the strict transform of a
  line meeting $\Cc$ in degree 1; the claim follows from $N_{1,1}=1$,
  where $N_{1,1}$ is the number of lines meeting a generic point and a
  generic line (see \cite{gottsche-pandharipande-1998}).
\end{proof}

\begin{lemma}[$M^3_k$ for $k=12,18,25$]
  Let $\Cc_2\subset\PP^3$ be a smooth rational curve of degree $2\leq
  d\leq 4$ which is disjoint from a line $\Cc_1\subset\PP^3$; let $M$
  be the blow-up of $\PP^3$ along $\Cc_1\sqcup\Cc_2$. Then
  $I_{L_0-F_1-F_2}(\pt)=d$.
\end{lemma}

\begin{proof}
  A curve of class $\beta=L_0-F_1-F_2$ is the strict transform of a
  line meeting $\Cc_1$ and $\Cc_2$ in one point each, or the connected
  union of the strict transform of a line which intersects $\Cc_1$ in
  one point and $\Cc_2$ in a scheme of degree $k>1$ with some
  exceptional fibers over $\Cc_2$. The class $\beta$ is enumerative,
  since reducible curves form at most a 1-dimensional family and
  irreducible ones form a 2-dimensional family.
  
  Consider the plane spanned by a generic point $x\in\PP^3$ and the
  line $\Cc_1$; it intersects $\Cc_2$ in a zero-dimensional scheme of
  degree $\deg\Cc_2$, and the intersection is transverse by the
  genericity assumptions on $x$.  This implies both that no reducible
  curves of class $\beta$ meet the strict transform $x'$ of $x$, and
  that exactly $\deg\Cc_2$ irreducible curves of class $\beta$ meet
  $x'$, which proves the lemma.
\end{proof}

\begin{lemma}\label{lemma:L0:3:10,15,20}
  $I_{L_0}(\ret,\pt)=1$ for $M^3_k$ with $k=10,15,20$.
\end{lemma}

\begin{proof}
  It is easy to see that the subspace of reducible curves has positive
  codimension in the moduli space $\Mzerozero(M,L_0)$; so 
  $$\dim\Mzerozero(M,L_0)=\dim\{\text{lines in the quadric}\}=3.$$
  The
  invariant $I_{L_0}(\ret,\pt)$ is thus enumerative, and counts the
  number of curves in the class $L_0$ which intersects the strict
  transforms both of a generic point $x'$ and a generic line $\ell'$
  in $Q^3$.
  
  Genericity assumptions imply that the unique $\ell$ meeting both
  $x'$ and $\ell'$ does not touch the blown-up locus; so the value of
  the invariant is $N_{1,1}=1$, as in lemma \ref{lemma:N11}.
\end{proof}

\endgroup


\begin{table}[h]
  \tablespace
  \centering
  \caption{Amount of geometric information needed in order to
    determine Quantum Cohomology. $L$ denotes a line, $C$ a conic, $T$
    a twisted cubic and $U$ a rational quartic.}
  \label{tab:geom-info}
  \bigskip
  \begin{tabular}{l|rrr}
    Threefold  & $N$ & $\dim A$ & $\deg A$\\\hline
    $\VSpace M^2_{22}=$ blow-up of $\PP^3$ along $U$         & 24 & 2 & 3 \\
    $\VSpace M^2_{27}=$ blow-up of $\PP^3$ along $T$         & 14 & 2 & 2 \\
    $\VSpace M^2_{30}=$ blow-up of $\PP^3$ along $C$         & 14 & 3 & 1 \\
    $\VSpace M^2_{33}=$ blow-up of $\PP^3$ along $L$         & 10 & 2 & 1 \\[2pt]\hline
    $\VSpace M^2_{21}=$ blow-up of $Q^3$ along $U$           & 24 & 3 & 5 \\
    $\VSpace M^2_{26}=$ blow-up of $Q^3$ along $T$           & 24 & 3 & 3 \\
    $\VSpace M^2_{29}=$ blow-up of $Q^3$ along $C$           & 14 & 3 & 2 \\[2pt]\hline
    $\VSpace M^3_{12}=$ blow-up of $\PP^3$ along $L\sqcup T$ & 81 & 3 & 9 \\
    $\VSpace M^3_{18}=$ blow-up of $\PP^3$ along $L\sqcup C$ & 81 & 3 & 6 \\
    $\VSpace M^3_{25}=$ blow-up of $\PP^3$ along $L\sqcup L$ & 52 & 3 & 3 \\[2pt]\hline
    $\VSpace M^3_{10}=$ blow-up of $Q^3$ along $C\sqcup C$   & 81 & 3 & 23\\
    $\VSpace M^3_{15}=$ blow-up of $Q^3$ along $L\sqcup C$   & 81 & 3 & 7 \\
    $\VSpace M^3_{20}=$ blow-up of $Q^3$ along $L\sqcup L$   & 81 & 3 & 5 
  \end{tabular}
\end{table}

\begin{table}[h]
  \tablespace
  \centering
  \caption{Relations in the small Quantum Cohomology rings of
    threefolds described in Theorem~\ref{thm:b2=2}.}
  \label{tab:rels-2}
  \bigskip
  \begin{tabular}{l|p{3.5cm}p{7.5cm}}
    $X$ & $f_i^C$ & $f_i^Q$
    \\\hline$
    \VSpace M^2_{21}
    $&$
    {E} H^{2}
    $&$
    -8 {E} {H} {{q}}_{0}+10 H^{2} {{q}}_{0}-28 {E} {{q}}_{0}^{2}-6 {H}
    {{q}}_{0}^{2}
    $\\&&$
    \ +8 {H} {{q}}_{0} {{q}}_{1}+16 {{q}}_{0}^{2} {{q}}_{1}%
    $\\[8pt]&$
    E^{2}-\tfrac52 {E} {H}+2 H^{2}
    $&$
    -2 {E} {{q}}_{0}+3 {H} {{q}}_{0}+{E} {{q}}_{1}+2 {{q}}_{0} {{q}}_{1}
    $\\[4pt]\hline$
    \VSpace M^2_{22}
    $&$
    {E} H^{2}
    $&$
    -6 {E} {H} {{q}}_{0}+10 H^{2} {{q}}_{0}-15 {E} {{q}}_{0}^{2}-6 {H}
    {{q}}_{0}^{2}
    $\\&&$
    \ -4 {E} {{q}}_{0} {{q}}_{1}+18 {H} {{q}}_{0} {{q}}_{1}+4 {{q}}_{0} {{q}}_{1}^{2}%
    $\\[8pt]&$
    E^{2}-\tfrac72 {E} {H}+4 H^{2}
    $&$
    -\tfrac52 {E} {{q}}_{0}+5 {H} {{q}}_{0}+{E} {{q}}_{1}+3 {{q}}_{0} {{q}}_{1}
    $\\[4pt]\hline$
    \VSpace M^2_{26}
    $&$
    {E} H^{2}
    $&$
    -4 {E} {H} {{q}}_{0}+\tfrac72 H^{2} {{q}}_{0}-6 {E} {{q}}_{0}^{2}-2 {H} {{q}}_{0}^{2}+6 {H} {{q}}_{0} {{q}}_{1}%
    $\\[4pt]&$
    E^{2}-\tfrac73 {E} {H}+\tfrac32 H^{2}
    $&$
    -\tfrac56 {E} {{q}}_{0}+\tfrac56 {H} {{q}}_{0}+{E} {{q}}_{1}+\tfrac12 {{q}}_{0} {{q}}_{1}
    $\\[4pt]\hline$
    \VSpace M^2_{27}
    $&$
    {E} H^{2}
    $&$
    -3Eq_0+8Hq_0+3q_0q_1
    $\\[4pt]&$
    E^2-\tfrac{10}3EH+3H^2
    $&$
    Eq_1+\tfrac13q_0
    $\\[4pt]\hline$
    \VSpace M^2_{29}
    $&$
    {E} H^{2}
    $&$
    4 {H} {{q}}_{0}
    $\\[4pt]&$
    E^{2}-2 {E} {H}+H^{2}
    $&$
    {E} {{q}}_{1}
    $\\[4pt]\hline$
    \VSpace M^2_{30}
    $&$
    {E} H^{2}
    $&$
    -2Eq_0+2Hq_0+2q_0q_1
    $\\[4pt]&$
    E^2-3EH+2H^2
    $&$
    Eq_1
    $\\[4pt]\hline$
    \VSpace M^2_{33}
    $&$
    {E} H^{2}
    $&$
    q_0
    $\\[4pt]&$
    E^2-2EH+H^2
    $&$
    Eq_1
    $
  \end{tabular}
\end{table}

\begin{table}[h]
  \tablespace
  \centering
  \caption{Relations in the small Quantum Cohomology rings of
    threefolds described in Theorem~\ref{thm:b2=3}.}
  \label{tab:rels-3}
  \bigskip
  \begin{tabular}{l|p{3.5cm}p{7.5cm}}
    $X$ & $f_i^C$ & $f_i^Q$
    \\\hline$
    \VSpace M^3_{10}
    $&$
    {{E}}_{1} {{E}}_{{2}}%
    $&$
    -2 {{E}}_{1} {{q}}_{0}-2 {{E}}_{{2}} {{q}}_{0}+4 {H} {{q}}_{0}+4 {{q}}_{0}^{2}%
    $\\[4pt]&$
    {{E}}_{1}^{2}-2 {{E}}_{1} {H}+H^{2}
    $&$
    {{E}}_{1} {{q}}_{1}+2 {{q}}_{0} {{q}}_{1}%
    $\\[4pt]&$
    {{E}}_{{2}}^{2}-2 {{E}}_{{2}} {H}+H^{2}    
    $&$
    {{E}}_{{2}} {{q}}_{{2}}+2 {{q}}_{0} {{q}}_{{2}}
    $\\[4pt]\hline$
    \VSpace M^3_{12}
    $&$
{{E}}_{1} {{E}}_{{2}}%
    $&$
-2 {{E}}_{1} {{q}}_{0}-4 {{E}}_{{2}} {{q}}_{0}+8 {H} {{q}}_{0}+3 {{q}}_{0} {{q}}_{{2}}%
    $\\[4pt]&$
{{E}}_{1}^{2}-2 {{E}}_{1} {H}+H^{2}%
    $&$
{{E}}_{1} {{q}}_{1}+{{q}}_{0} {{q}}_{1}%
    $\\[4pt]&$
{{E}}_{{2}}^{2}-\tfrac{10}3 {{E}}_{{2}} {H}+3 H^{2}
    $&$
-\tfrac13 {{E}}_{1} {{q}}_{0}-\tfrac23 {{E}}_{{2}} {{q}}_{0}+\tfrac43 {H} {{q}}_{0}+\tfrac13 {{q}}_{0} {{q}}_{1}+{{E}}_{{2}} {{q}}_{{2}}+2 {{q}}_{0} {{q}}_{{2}}
    $\\[4pt]\hline$
    \VSpace M^3_{15}
    $&$
{{E}}_{1} {{E}}_{{2}}%
    $&$
-2 {{E}}_{1} {{q}}_{0}-{{E}}_{{2}} {{q}}_{0}+2 {H} {{q}}_{0}%
    $\\[4pt]&$
{{E}}_{1}^{2}-{{E}}_{1} {H}+\tfrac12 H^{2}%
    $&$
-{{E}}_{1} {{q}}_{0}-\tfrac12 {{E}}_{{2}} {{q}}_{0}+{H} {{q}}_{0}+{{E}}_{1} {{q}}_{1}+{{q}}_{0} {{q}}_{1}+\tfrac12 {{q}}_{0} {{q}}_{{2}}%
    $\\[4pt]&$
{{E}}_{{2}}^{2}-2 {{E}}_{{2}} {H}+H^{2}
    $&$
{{E}}_{{2}} {{q}}_{{2}}+{{q}}_{0} {{q}}_{{2}}
    $\\[4pt]\hline$
    \VSpace M^3_{18}
    $&$
{{E}}_{1} {{E}}_{{2}}%
    $&$
-2 {{E}}_{{2}} {{q}}_{0}+2 {H} {{q}}_{0}+2 {{q}}_{0} {{q}}_{{2}}%
    $\\[4pt]&$
{{E}}_{1}^{2}-2 {{E}}_{1} {H}+H^{2}%
    $&$
{{E}}_{1} {{q}}_{1}%
    $\\[4pt]&$
{{E}}_{{2}}^{2}-3 {{E}}_{{2}} {H}+2 H^{2}
    $&$
{{E}}_{{2}} {{q}}_{{2}}
    $\\[4pt]\hline$
    \VSpace M^3_{20}
    $&$
{{E}}_{1} {{E}}_{{2}}%
    $&$
-{{E}}_{1} {{q}}_{0}-{{E}}_{{2}} {{q}}_{0}+{H} {{q}}_{0}%
    $\\[4pt]&$
{{E}}_{1}^{2}-{{E}}_{1} {H}+\tfrac12 H^{2}%
    $&$
-\tfrac12 {{E}}_{1} {{q}}_{0}-\tfrac12 {{E}}_{{2}} {{q}}_{0}+\tfrac12 {H} {{q}}_{0}+{{E}}_{1} {{q}}_{1}+\tfrac12 {{q}}_{0} {{q}}_{1}+\tfrac12 {{q}}_{0} {{q}}_{{2}}%
    $\\[4pt]&$
{{E}}_{{2}}^{2}-{{E}}_{{2}} {H}+\tfrac12 H^{2}
    $&$
-\tfrac12 {{E}}_{1} {{q}}_{0}-\tfrac12 {{E}}_{{2}} {{q}}_{0}+\tfrac12 {H} {{q}}_{0}+\tfrac12 {{q}}_{0} {{q}}_{1}+{{E}}_{{2}} {{q}}_{{2}}+\tfrac12 {{q}}_{0} {{q}}_{{2}}
    $\\[4pt]\hline$
    \VSpace M^3_{25}
    $&$
{{E}}_{1} {{E}}_{{2}}%
    $&$
{{q}}_{0}%
    $\\[4pt]&$
{{E}}_{1}^{2}-2 {{E}}_{1} {H}+H^{2}%
    $&$
{{E}}_{1} {{q}}_{1}%
    $\\[4pt]&$
{{E}}_{{2}}^{2}-2 {{E}}_{{2}} {H}+H^{2}%
    $&$
{{E}}_{{2}} {{q}}_{{2}}%
    $
  \end{tabular}
\end{table}


\bibliographystyle{amsplain} 

\bibliography{qh}

\def\cprime{$'$}
\providecommand{\bysame}{\leavevmode\hbox to3em{\hrulefill}\thinspace}
\providecommand{\MR}{\relax\ifhmode\unskip\space\fi MR }
\providecommand{\MRhref}[2]{%
  \href{http://www.ams.org/mathscinet-getitem?mr=#1}{#2}
}
\providecommand{\href}[2]{#2}
\begin{thebibliography}{10}

\bibitem{ancona-maggesi-2002}
Vincenzo Ancona and Marco Maggesi, \emph{{Quantum Cohomology} of some {Fano}
  threefolds}, to appear in {\em Advances in Geometry}, 2002.

\bibitem{ancona-ottaviani-1989}
Vincenzo Ancona and Giorgio Ottaviani, \emph{An introduction to the derived
  categories and the theorem of {B}eilinson}, Atti Accad. Peloritana
  Pericolanti Cl. Sci. Fis. Mat. Natur. \textbf{67} (1989), 99--110 (1991).
  \MR{92g:14013}

\bibitem{ancona-ottaviani-1991}
\bysame, \emph{Canonical resolutions of sheaves on {S}chubert and {B}rieskorn
  varieties}, Complex analysis (Wuppertal, 1991), Aspects Math., no. E17,
  Vieweg, Braunschweig, 1991, pp.~14--19.

\bibitem{algheme}
Giovanni Baldini, Gianni Ciolli, and Marco Maggesi, \emph{Algheme, a
  commutative algebra package for the {S}cheme language}, Freely available at
  {\tt http://sourceforge.net/projects/algheme/}.

\bibitem{batyrev-1981}
V.~V. Batyrev, \emph{Toric {F}ano threefolds}, Izv. Akad. Nauk SSSR Ser. Mat.
  \textbf{45} (1981), no.~4, 704--717, 927. \MR{82m:14022}

\bibitem{batyrev-1993}
Victor~V. Batyrev, \emph{Quantum cohomology rings of toric manifolds},
  Ast\'erisque (1993), no.~218, 9--34, Journ\'ees de G\'eom\'etrie Alg\'ebrique
  d'Orsay (Orsay, 1992). \MR{95b:32034}

\bibitem{bayer-2004}
Arend Bayer, \emph{{S}emisimple {Q}uantum {C}ohomology and {B}low-ups},
  Preprint {\tt arXiv:math.AG/0403260}, 2004.

\bibitem{bayer-manin-2001}
Arend Bayer and Yuri~I. Manin, \emph{{(Semi)simple} exercises in {Quantum
  Cohomology}}, Preprint {\tt arXiv:math.AG/0103164}, 2001.

\bibitem{beauville-1995}
A.~Beauville, \emph{Quantum cohomology of complete intersections}, Mat. Fiz.
  Anal. Geom. \textbf{2} (1995), no.~3-4, 384--398.

\bibitem{behrend-1997}
K.~Behrend, \emph{Gromov-{W}itten invariants in algebraic geometry}, Invent.
  Math. \textbf{127} (1997), no.~3, 601--617. \MR{98i:14015}

\bibitem{behrend-fantechi-1997}
K.~Behrend and B.~Fantechi, \emph{The intrinsic normal cone}, Invent. Math.
  \textbf{128} (1997), no.~1, 45--88. \MR{98e:14022}

\bibitem{beilinson-1978}
A.~A. Beilinson, \emph{Coherent sheaves on {${\bf P}\sp{n}$} and problems in
  linear algebra}, Funktsional. Anal. i Prilozhen. \textbf{12} (1978), no.~3,
  68--69. \MR{80c:14010b}

\bibitem{fano3fold-software}
Gianni Ciolli, \emph{{\tt fano3folds}, a software for the systematic
  computation of {A}ssociativity relations for {F}ano threefolds}, can be
  downloaded freely at {\tt
  http://www.math.unifi.it/{\textasciitilde}ciolli/fano3folds}, 2003.

\bibitem{ciolli-phd-thesis}
\bysame, \emph{On the {Q}uantum {C}ohomology of some {F}ano threefolds and a
  conjecture of {D}ubrovin}, Tesi di Dottorato, 2003.

\bibitem{costa-miro-roig-2000}
L.~Costa and R.~M. Mir\`o-Roig, \emph{Quantum cohomology of projective bundles
  over {$\PP^{n_1}\times\dots\times\PP^{n_s}$}}, International J. of Math.
  \textbf{11} (2000), no.~6, 761--797.

\bibitem{dubrovin-1998}
Boris Dubrovin, \emph{Geometry and analytic theory of {F}robenius manifolds},
  Proceedings of the International Congress of Mathematicians, Vol. II, Doc.
  Math. 1998, no. Extra Vol. II, Berlin, 1998, pp.~315--326.

\bibitem{faenzi-2003-v22}
Daniele Faenzi, \emph{Bundles over {F}ano {T}hreefolds of type {$V_{22}$}},
  quaderno del {D}ipartimento di {M}atematica ``{U}lisse {D}ini'', 2003.

\bibitem{fulton-pandharipande-1997}
W.~Fulton and R.~Pandharipande, \emph{Notes on stable maps and quantum
  cohomology}, Algebraic geometry---Santa Cruz 1995, Proc. Sympos. Pure Math.,
  vol.~62, Amer. Math. Soc., Providence, RI, 1997, pp.~45--96. \MR{98m:14025}

\bibitem{fulton-1984}
William Fulton, \emph{Intersection theory}, Ergebnisse der Mathematik und ihrer
  Grenzgebiete (3) [Results in Mathematics and Related Areas (3)], vol.~2,
  Springer-Verlag, Berlin, 1984. \MR{85k:14004}

\bibitem{givental-1998}
Alexander Givental, \emph{A mirror theorem for toric complete intersections},
  Topological field theory, primitive forms and related topics (Kyoto, 1996),
  Progr. Math., vol. 160, Birkh\"auser Boston, Boston, MA, 1998, pp.~141--175.
  \MR{2000a:14063}

\bibitem{gottsche-pandharipande-1998}
L.~G\"ottsche and R.~Pandharipande, \emph{The {Quantum Cohomology} of blow-ups
  of {$\PP^2$} and enumerative geometry}, J. Differential Geom. \textbf{48}
  (1998), no.~1, 61--90.

\bibitem{hu-2000}
J.~Hu, \emph{Gromov-{W}itten invariants of blow-ups along points and curves},
  Math. Z. \textbf{233} (2000), no.~4, 709--739. \MR{2001c:53115}

\bibitem{iskovskih-1977}
V.~A. Iskovskih, \emph{Fano threefolds. {I}}, Izv. Akad. Nauk SSSR Ser. Mat.
  \textbf{41} (1977), no.~3, 516--562, 717. \MR{80c:14023a}

\bibitem{iskovskih-1978}
\bysame, \emph{Fano threefolds. {II}}, Izv. Akad. Nauk SSSR Ser. Mat.
  \textbf{42} (1978), no.~3, 506--549. \MR{80c:14023b}

\bibitem{kontsevich-manin-1994}
M.~Kontsevich and Yuri~I. Manin, \emph{Gromov-{W}itten classes, quantum
  cohomology, and enumerative geometry}, Comm. Math. Phys. \textbf{164} (1994),
  no.~3, 525--562. \MR{95i:14049}

\bibitem{kuznetsov-1996}
A.~G. Kuznetsov, \emph{An exception set of vector bundles on the varieties
  {$V\sb {22}$}}, Vestnik Moskov. Univ. Ser. I Mat. Mekh. (1996), no.~3,
  41--44, 92. \MR{97m:14040}

\bibitem{li-tian-1998}
Jun Li and Gang Tian, \emph{Virtual moduli cycles and {G}romov-{W}itten
  invariants of algebraic varieties}, J. Amer. Math. Soc. \textbf{11} (1998),
  no.~1, 119--174. \MR{99d:14011}

\bibitem{mori-mukai-1982}
Shigefumi Mori and Shigeru Mukai, \emph{Classification of {F}ano {$3$}-folds
  with {$B\sb{2}\geq 2$}}, Manuscripta Math. \textbf{36} (1981/82), no.~2,
  147--162. \MR{83f:14032}

\bibitem{mori-mukai-2003}
\bysame, \emph{Erratum to ``classification of {F}ano {$3$}-folds with
  {$B\sb{2}\geq 2$}''}, Manuscripta Math. \textbf{110} (2003), 407.

\bibitem{orlov-1991}
D.~O. Orlov, \emph{Exceptional set of vector bundles on the variety {$V\sb
  5$}}, Vestnik Moskov. Univ. Ser. I Mat. Mekh. (1991), no.~5, 69--71.
  \MR{95f:14080}

\bibitem{orlov-1993}
D.~O. Orlov, \emph{Projective bundles, monoidal transformations, and derived
  categories of coherent sheaves}, Russian Acad. Sci. Izv. Math. \textbf{41}
  (1993), no.~1, 133--141.

\bibitem{qin-ruan-1998}
Z.~Qin and Y.~Ruan, \emph{Quantum cohomology of projective bundles over
  {$\PP^n$}}, Transactions of the Am. Math. Soc. \textbf{350} (1998), no.~9,
  3615--3638.

\bibitem{ruan-tian-1995}
Yongbin Ruan and Gang Tian, \emph{A mathematical theory of quantum cohomology},
  J. Differential Geom. \textbf{42} (1995), no.~2, 259--367. \MR{96m:58033}

\bibitem{spielberg-1999}
Holger Spielberg, \emph{The {G}romov-{W}itten invariants of symplectic toric
  manifolds, and their quantum cohomology ring}, C. R. Acad. Sci. Paris S\'er.
  I Math. \textbf{329} (1999), no.~8, 699--704. \MR{2000j:14087}

\bibitem{szurek-wisniewski-1990}
Micha{\l} Szurek and Jaros{\l}aw~A. Wi{\'s}niewski, \emph{Fano bundles of rank
  {$2$} on surfaces}, Compositio Math. \textbf{76} (1990), no.~1-2, 295--305,
  Algebraic geometry (Berlin, 1988). \MR{92e:14037}

\bibitem{watanabe-watanabe-1982}
Keiichi Watanabe and Masayuki Watanabe, \emph{The classification of {F}ano
  {$3$}-folds with torus embeddings}, Tokyo J. Math. \textbf{5} (1982), no.~1,
  37--48. \MR{83m:14029}

\end{thebibliography}

\end{document}